\documentclass[8pt]{article}

\usepackage{amsmath,amssymb,amsfonts,amsthm,amscd,shadow,fancybox,epic}
\usepackage{pstricks,pst-node,pst-plot}

\theoremstyle{plain}

\newtheorem{lem}[subsection]{Lemma}
\newtheorem{prop}[subsection]{Proposition}

\newtheorem{ex}[subsection]{Example}

\theoremstyle{definition}
\newtheorem{rem}[subsection]{Remark}

\newcommand{\pa}{\zeta}
\newcommand{\pb}{\eta}
\newcommand{\spa}{\lambda}
\newcommand{\spb}{\mu}

\newcommand{\papb}{_{\pa,\pb}}
\newcommand{\spaspb}{_{\spa,\spb}}

\newcommand{\gl}{\lambda}
\newcommand{\gm}{\mu}
\newcommand{\gp}{\pi}
\newcommand{\gt}{\theta}
\newcommand{\ga}{\alpha}
\newcommand{\gb}{\beta}

\newcommand{\gr}{\rho}
\newcommand{\gs}{\sigma}

\newcommand{\ft}{\mathfrak{t}}
\newcommand{\fR}{\mathfrak{R}}
\newcommand{\bC}{\mathbb{C}}

\newcommand{\cF}{\mathcal{F}}
\newcommand{\cD}{\mathcal{D}}

\newcommand{\supp}{\mathop{\rm Supp}\nolimits}

\newcommand{\Span}{\mathop{\rm Span}\nolimits}
\newcommand{\Char}{\mathop{\rm Char}\nolimits}

\newcommand{\diag}{\mathop{\rm diag}\nolimits}

\newcommand{\ssvt}{\mathop{\rm SSV\widetilde{T}}\nolimits}
\newcommand{\ssyt}{\mathop{\rm SSY\widetilde{T}}\nolimits}
\newcommand{\svt}{\mathop{\rm SV\widetilde{T}}\nolimits}
\newcommand{\ossyt}{\mathop{\rm SSY\overline{T}}\nolimits}

\newcommand{\ol}{\overline}

\newcommand{\sK}{_{\text{\tiny{K}}}}
\newcommand{\sH}{_{\text{\tiny{H}}}}

\newcommand{\cO}{\mathcal{O}}
\newcommand{\ab}{_{\ga,\gb}}
\newcommand{\lm}{_{\gl,\gm}}

\newcommand{\ocD}{\overline{\mathcal{D}}}
\newcommand{\ocF}{\overline{\mathcal{F}}}
\newcommand{\wcD}{\widetilde{\mathcal{D}}}
\newcommand{\wcF}{\widetilde{\mathcal{F}}}
\newcommand{\wt}{\widetilde}

\title{Schubert Classes in the
Equivariant K-Theory and Equivariant Cohomology of the\\
Lagrangian Grassmannian}

\author{Victor Kreiman}

\date{}

\begin{document}
\maketitle

\begin{abstract}
We give positive formulas for the restriction of a Schubert Class
to a $T$-fixed point in the equivariant K-theory and equivariant
cohomology of the Lagrangian Grassmannian.  Our formulas rely on a
result of Ghorpade-Raghavan, which gives an equivariant Gr\"obner
degeneration of a Schubert variety in the neighborhood of a
$T$-fixed point of the Lagrangian Grassmannian.
\end{abstract}

\section{Introduction}\label{s.intro}

Let $J$ be the antidiagonal $2n\times 2n$ matrix whose top $n$
antidiagonal entries are 1's and whose bottom $n$ antidiagonal
entries are -1's. Then $J$ defines a nondegenerate skew-symmetric
inner product on $\bC^{2n}$ by $\langle v,w\rangle=v^t J w$,
$v,w\in\bC^{2n}$. The Lagrangian Grassmannian $LGr_n$ is defined
as the set of all $n$-dimensional complex subspaces $V$ of
$\bC^{2n}$ which are isotropic under this inner product, i.e.,
such that for every $v,w\in V$, $\langle v,w\rangle=0$. The
symplectic group $G=Sp_{2n}(\bC)$ consists of the invertible
$2n\times 2n$ complex matrices which preserve this inner product.
Let $T$ and $B$ denote the diagonal and upper triangular matrices
of $G$ respectively.  The natural action of $G$ on $LGr_{n}$ is
transitive and has a unique $B$-fixed point $e_{id}$. Thus $LGr_n$
can be identified with $G/P_n$, where $P_n\supset B$ is the
stabilizer of $e_{id}$. Let $W$ denote the Weyl group of $G$ with
respect to $T$ ($=N_G(T)/T$) and $W_{P_n}$ the Weyl group of
$P_n$. For the $G$-action on $LGr_n$, the $T$-fixed points are
precisely the cosets $e_\gb:=\gb P_n$, $\gb\in W/W_{P_n}$.

Let $B^-$ denote the lower triangular matrices in $G$.   For
$\ga\in W/W_{P_n}$, the (opposite) Schubert variety $X_\ga$ is the
Zariski closure of $B^- e_\ga$ in $LGr_n$. The Schubert variety
$X_\ga$ defines classes $[X_\ga]\sK$ in $K_T^*(LGr_n)$, the
$T$-equivariant K-theory of $LGr_n$, and $[X_\ga]\sH$ in
$H_T^*(LGr_n)$, the $T$-equivariant cohomology of $LGr_n$.

The $T$-equivariant embedding $e_\gb\stackrel{i}{\to} LGr_n$
induces restriction homomorphisms:
\begin{equation*}
K_T^*(LGr_{n})\stackrel{i^*_K}{\rightarrow} K_T^*(e_\gb)
\qquad\text{and}\qquad
H_T^*(LGr_{n})\stackrel{i^*_H}{\rightarrow} H_T^*(e_\gb).
\end{equation*}
The image of an element $C$ of $K_T^*(LGr_{n})$ or
$H_T^*(LGr_{n})$ under restriction to $e_\gb$ is denoted by
$C|_{e_\gb}$. The restrictions $C|_{e_\gb}$, evaluated at all
$\gb\in W/W_{P_n}$, determine $C$ uniquely. In this paper we
obtain combinatorial formulas for $[X_\ga]\sK|_{e_\gb}$ and
$[X_\ga]\sH|_{e_\gb}$. Our formula for $[X_\ga]\sK|_{e_\gb}$ is
positive in the sense of \cite[Conjecture 5.1]{Gr-Ra}, and our
formula for $[X_\ga]\sH|_{e_\gb}$ is positive in the sense of
\cite{Gr}. A positive formula for $[X_\ga]\sK|_{e_\gb}$ also
appears in \cite{Wi}, and positive formulas for
$[X_\ga]\sH|_{e_\gb}$ appear in \cite{Bi} and \cite{Ik}.

The proof of our formulas relies on a result of Ghorpade-Raghavan
\cite{Gh-Ra}, which gives an explicit equivariant Gr\"obner
degeneration of an open neighborhood of $X_\ga$ centered at
$e_\gb$ to a reduced union of coordinate spaces. The outline of
our proof is virtually the same as that of \cite{Kr}, which
derives a similar result as here, but for Schubert varieties in
the ordinary Grassmannian. In addition, many of the lemmas of
\cite{Kr} and their proofs carry over with little or no
modification.

Our formulas for $[X_\ga]\sK|_{e_\gb}$ and $[X_\ga]\sH|_{e_\gb}$
are expressed in terms of `semistandard set-valued {\em shifted}
tableaux'. These objects take the place of the `semistandard
set-valued tableaux' in \cite{Kr}. Semistandard set-valued
tableaux were introduced by Buch \cite{Bu}, and also appear in
\cite{Kn-Mi-Yo1}, \cite{Kn-Mi-Yo2}. The formula for
$[X_\ga]\sH|_{e_\gb}$ can also be expressed in terms of `subsets
of {\em shifted} diagrams', which we introduce in Section
\ref{s.three_equiv_models}. These objects take the place of the
`subsets of Young diagrams' in \cite{Kr}.
It has come to our attention that Ikeda-Naruse have independently
discovered subsets of Young diagrams and subsets of shifted
diagrams and used them to express formulas for restrictions of
Schubert classes to $T$-fixed points in the equivariant cohomology
of the ordinary and Lagrangian Grassmannians respectively.
\vspace{1em}

\section{Semistandard Set-Valued Shifted Tableaux}\label{s.ssvt}

For $k\in \{1,\ldots,2n\}$, define $\ol{k}=2n+1-k$. Let $I_n$
denote the set of all $n$ element subsets
$\ga=\{\ga(1),\ldots,\ga(n)\}$ of $\{1,\ldots,2n\}$ such that for
each $k\in \{1,\ldots,2n\}$, exactly one of $k$ or $\ol{k}$ is in
$\ga$. We always assume the entries of such a subset are listed in
increasing order. For $\ga\in I_n$, define ${\ga'}\in I_{n}$ by
${\ga'}=\{1,\ldots,2n\}\setminus
\ga=\{\ol{\ga(n)},\ldots,\ol{\ga(1)}\}$. The map which takes
$\{\ga(1),\ldots,\ga(n)\}\in I_n$ to the permutation
$(\ga(1),\ldots,\ga(n),\ol{\ga(n)},\ldots,\ol{\ga(1)})\in W$
identifies $I_n$ with the set of minimal length coset
representatives for $W/W_{P_n}$. We shall use $I_n$ rather than
$W/W_{P_n}$ to index the Schubert varieties and $T$-fixed points
of $LGr_n$. Fix $\ga,\gb\in I_n$ for the remainder of this paper.

A \textbf{partition} is an ordered list of nonnegative integers
$\gl=(\gl_1,\ldots,\gl_m)$, $\gl_1\geq\cdots\geq\gl_m$. Two
partitions are identified if one can be obtained from the other by
adding zeros. The \textbf{transpose} of $\gl$ is the partition
$\gl^t=(\gl^t_1,\ldots,\gl^t_p)$, where
$\gl^t_j=\#\{i\in\{1,\ldots,m\}\mid\gl_i\geq j\}$, $j=1,\ldots,p$.
The partition $\gl$ is said to be \textbf{symmetric} if
$\gl^t=\gl$. It is said to be \textbf{strict} if $\gl_i=\gl_{i+1}$
implies $\gl_i=0$, $i=1,\ldots,m$. We denote by $L_n$ (resp.
$M_n$) the set of all symmetric (resp. strict) partitions $\gl$
with $\gl_1\leq n$.

The map $\gp$ from finite subsets of the positive integers to
partitions, given by $\gp:\{\gamma(1),\ldots,\gamma(k)\}\mapsto
(\gamma(k)-k,\ldots,\gamma(1)-1)$, where
$\gamma(1)<\cdots<\gamma(k)$, restricts to a bijection from $I_n$
to $L_n$. The map $\gr$ from partitions to strict partitions,
given by $\gr:(\gl_1,\ldots,\gl_k)\mapsto
(\gl_1,\gl_2-1,\ldots,\gl_l-l+1)$, where $l$ is maximal such that
$\gl_l-l+1\geq 0$, restricts to a bijection from $L_n$ to $M_n$.
We denote the composition $\gr\circ\gp:I_n\to M_n$ by $\gs$.  If
$\gl=\gs(\ga)$, then the \textbf{length} of $\ga$, denoted
$l(\ga)$, is $\gl_1+\cdots+\gl_n$.

A \textbf{Young diagram} is a collection of boxes arranged into a
top and left justified array. A Young diagram is said to be
\textbf{symmetric} if the length of the $i$-th row equals the
length of the $i$-th column for all $i$. To any partition $\gl$ we
associate the Young diagram $D_\gl$ whose $i$-th row has length
$\gl_i$. The $j$-th column of $D_\gl$ has length $\gl^t_j$. Thus
$\gl$ is symmetric if and only if $D_\gl$ is symmetric, and $L_n$
can be identified with the set of all symmetric Young diagrams
whose first rows have length $\leq n$.

A \textbf{shifted diagram} is a top-justified array of boxes whose
left side forms a descending staircase, i.e., the leftmost box of
any row is one column to the right of the leftmost box of the row
above it. The \textbf{length} of a row of a shifted diagram is the
number of boxes it contains. To a strict partition $\gl$ we
associate the shifted diagram $\wt{D}_\gl$ whose $i$-th row has
length $\gl_i$. We call $\gl$ the \textbf{shape} of $\wt{D}_\gl$.
 One sees that $M_n$ can be identified with the set of all shifted
diagrams whose first rows have length $\leq n$.

\psset{unit=.5cm}
\begin{equation*}
\begin{array}{c@{\hspace{3em}}c}
\pspicture*(0,0)(5,5)
\psset{linewidth=.04}
\psline(0,5)(5,5)
\psline(0,4)(5,4)
\psline(0,3)(3,3)
\psline(0,2)(2,2)
\psline(0,1)(1,1)
\psline(0,0)(1,0)
\psline(0,0)(0,5)
\psline(1,0)(1,5)
\psline(2,2)(2,5)
\psline(3,3)(3,5)
\psline(4,4)(4,5)
\psline(5,4)(5,5)
\endpspicture
&
\pspicture*(0,0)(10,5)
\psset{linewidth=.04}
\psline(0,5)(10,5)
\psline(0,4)(10,4)
\psline(1,3)(9,3)
\psline(2,2)(7,2)
\psline(3,1)(5,1)
\psline(4,0)(5,0)
\psline(0,4)(0,5)
\psline(1,3)(1,5)
\psline(2,2)(2,5)
\psline(3,1)(3,5)
\psline(4,0)(4,5)
\psline(5,0)(5,5)
\psline(6,2)(6,5)
\psline(7,2)(7,5)
\psline(8,3)(8,5)
\psline(9,3)(9,5)
\psline(10,4)(10,5)
\endpspicture
\\[.5em]
%
\hbox{(a) A symmetric Young diagram}
&
\hbox{(b) A shifted diagram}
\end{array}
\end{equation*}

The bijection $\gr:L_n\to M_n$ can be viewed in terms of
associated partitions. Let $\gl$ be a symmetric partition.  If we
remove all boxes of $D_\gl$ which lie below the main diagonal,
then we obtain $\wt{D}_{\gr(\gl)}$: \vspace{1em}

\begin{center}
\pspicture(0,0)(15,5)
\psset{linewidth=.04}
\psline(0,5)(5,5)
\psline(0,4)(5,4)
\psline(0,3)(3,3)
\psline(0,2)(2,2)
\psline(0,1)(1,1)
\psline(0,0)(1,0)
\psline(0,0)(0,5)
\psline(1,0)(1,5)
\psline(2,2)(2,5)
\psline(3,3)(3,5)
\psline(4,4)(4,5)
\psline(5,4)(5,5)
\psline[linestyle=dashed,dash=3pt
2pt,linewidth=.1](0,4)(1,4)(1,3)(2,3)
\rput{0}(7.5,3){$\longrightarrow$}
\psline(10,5)(15,5)
\psline(10,4)(15,4)
\psline(11,3)(13,3)
\psline(10,4)(10,5)
\psline(11,3)(11,5)
\psline(12,3)(12,5)
\psline(13,3)(13,5)
\psline(14,4)(14,5)
\psline(15,4)(15,5)
\endpspicture
\end{center}
\vspace{.1em}
\begin{center}
Figure 1: The map $\gr$
\end{center}

A \textbf{set-valued shifted tableau} $S$ is an assignment of a
nonempty set of positive integers to each box of a shifted
diagram. The entries of $S$ are the positive integers in the
boxes. If a positive integer occurs in more than one box of $S$,
then we consider the separate occurrences to be distinct entries.
If $x$ is an entry of $S$, then we define $r(x)$ and $c(x)$ to be
the row and column numbers of the box containing $x$ (where the
top row is considered the first row and the leftmost column is
considered the first column), and we define $z(x)$ to be
$x+c(x)-r(x)$. We say that $S$ is a \textbf{Young shifted tableau}
if each box contains a single entry.

A set-valued shifted tableau is said to be \textbf{semistandard}
if all entries of any box $B$ are less than or equal to all
entries of the box to the right of $B$ and strictly less than all
entries of the box below $B$. \vspace{1em}

\psset{unit=1.3cm}
\begin{center}
\pspicture(0,0)(6,3)
\psset{linewidth=.02}
\psline(0,3)(6,3)
\psline(0,2)(6,2)
\psline(1,1)(5,1)
\psline(2,0)(3,0)
\psline(0,2)(0,3)
\psline(1,1)(1,3)
\psline(2,0)(2,3)
\psline(3,0)(3,3)
\psline(4,1)(4,3)
\psline(5,1)(5,3)
\psline(6,2)(6,3)
\rput{0}(.5,2.5){$1$}
\rput{0}(1.5,2.5){$2,3$}
\rput{0}(2.5,2.5){$3$}
\rput{0}(3.5,2.5){$3$}
\rput{0}(4.5,2.5){$4,6,7$}
\rput{0}(5.5,2.5){$7,9$}
\rput{0}(1.5,1.5){$4$}
\rput{0}(2.5,1.5){$4,6$}
\rput{0}(3.5,1.5){$6,7,8$}
\rput{0}(4.5,1.5){$9,11$}
\rput{0}(2.5,0.5){$8,10$}
\endpspicture
\end{center}

\begin{center}
Figure 2: A semistandard set-valued shifted tableau
\end{center}
\vspace{1em}

\noindent If $\spb=(\spb_1,\ldots,\spb_h)$ is any strict
partition, then a set-valued shifted tableau $S$ is said to be
\textbf{on $\spb$} if, for every entry $x$ of $S$, $x\leq h$ and
\begin{equation}\label{e.semi_stand_on_mu}
z(x)\leq \spb_x+x-1.
\end{equation}

\begin{ex}\label{ex.ss_shifted_tabl} Let $\spa=(2,1)$, $\spb=(5,3,2)$. The following
list gives all semistandard set-valued shifted tableaux on $\spb$
of shape $\spa$:
\vspace{1em}
\psset{unit=1cm}
\begin{equation*}
\pspicture(0,0)(2,2)
\psset{linewidth=.024}
\psline(0,2)(2,2)
\psline(0,1)(2,1)
\psline(1,0)(2,0)
\psline(0,1)(0,2)
\psline(1,0)(1,2)
\psline(2,0)(2,2)
\rput{0}(.5,1.5){$1$}
\rput{0}(1.5,1.5){$1$}
\rput{0}(1.5,0.5){$2$}
\endpspicture
\qquad\qquad
\pspicture(0,0)(2,2)
\psset{linewidth=.024}
\psline(0,2)(2,2)
\psline(0,1)(2,1)
\psline(1,0)(2,0)
\psline(0,1)(0,2)
\psline(1,0)(1,2)
\psline(2,0)(2,2)
\rput{0}(.5,1.5){$1$}
\rput{0}(1.5,1.5){$1$}
\rput{0}(1.5,0.5){$3$}
\endpspicture
\qquad\qquad
\pspicture(0,0)(2,2)
\psset{linewidth=.024}
\psline(0,2)(2,2)
\psline(0,1)(2,1)
\psline(1,0)(2,0)
\psline(0,1)(0,2)
\psline(1,0)(1,2)
\psline(2,0)(2,2)
\rput{0}(.5,1.5){$1$}
\rput{0}(1.5,1.5){$2$}
\rput{0}(1.5,0.5){$3$}
\endpspicture
\qquad\qquad
\pspicture(0,0)(2,2)
\psset{linewidth=.024}
\psline(0,2)(2,2)
\psline(0,1)(2,1)
\psline(1,0)(2,0)
\psline(0,1)(0,2)
\psline(1,0)(1,2)
\psline(2,0)(2,2)
\rput{0}(.5,1.5){$2$}
\rput{0}(1.5,1.5){$2$}
\rput{0}(1.5,0.5){$3$}
\endpspicture
\end{equation*}
\vspace{.2pt}
\begin{equation*}
\quad
\pspicture(0,0)(2,2)
\psset{linewidth=.024}
\psline(0,2)(2,2)
\psline(0,1)(2,1)
\psline(1,0)(2,0)
\psline(0,1)(0,2)
\psline(1,0)(1,2)
\psline(2,0)(2,2)
\rput{0}(.5,1.5){$1$}
\rput{0}(1.5,1.5){$1$}
\rput{0}(1.5,0.5){$2,3$}
\endpspicture
\qquad\qquad
\pspicture(0,0)(2,2)
\psset{linewidth=.024}
\psline(0,2)(2,2)
\psline(0,1)(2,1)
\psline(1,0)(2,0)
\psline(0,1)(0,2)
\psline(1,0)(1,2)
\psline(2,0)(2,2)
\rput{0}(.5,1.5){$1$}
\rput{0}(1.5,1.5){$1,2$}
\rput{0}(1.5,0.5){$3$}
\endpspicture
\qquad\qquad
\pspicture(0,0)(2,2)
\psset{linewidth=.024}
\psline(0,2)(2,2)
\psline(0,1)(2,1)
\psline(1,0)(2,0)
\psline(0,1)(0,2)
\psline(1,0)(1,2)
\psline(2,0)(2,2)
\rput{0}(.5,1.5){$1,2$}
\rput{0}(1.5,1.5){$2$}
\rput{0}(1.5,0.5){$3$}
\endpspicture
\end{equation*}
\end{ex}
\vspace{.5em}

\noindent Denote the set of semistandard set-valued shifted
tableaux on $\spb$ of shape $\spa$ by $\ssvt\spaspb$ and the set of
semistandard Young shifted tableaux on $\spb$ of shape $\spa$ by
$\ssyt\spaspb$.

\section{Results}\label{s.results}

Let $\ft$ denote the Lie algebra of $T$ and $R(T)$ the
representation ring of $T$. We have that
\begin{align*}
&T=\{\diag(s_1,\ldots,s_{n},s_n^{-1},\ldots,s_1^{-1})\mid
s_k\in\bC^*\}\\
&\ft=\{\diag(s_1,\ldots,s_{n},-s_n,\ldots,-s_1)\mid
s_k\in\bC\}\\[.5em]
&K_T^*(e_\gb)\cong R(T)=\bC[t_1^{\pm 1},\ldots,t_n^{\pm 1}]
\\
&H_T^*(e_\gb)\cong \bC[\ft^*] =\bC[t_1,\ldots,t_n]
\end{align*}
For $k=1,\ldots,n$, define $t_{\ol{k}}\in K_T^*(e_\gb)$ to be
$t_k^{-1}$ and $t_{\ol{k}}\in H_T^*(e_\gb)$ to be $-t_k$.

\begin{prop}\label{p.rest_formula_kthry}
Let $\spa=\gs(\ga)$, $\spb=\gs(\gb)$. Then

\noindent (i) $\displaystyle
[X_\ga]\sK|_{e_\gb}=(-1)^{l(\ga)}\sum\limits_{S\in\ssvt\spaspb} \,
\prod\limits_{x\in
S}\left(\frac{1}{t_{{\gb'}(x)}t_{{\gb'}(z(x))}}-1\right)$.
\vspace{1em}

\noindent (ii) $\displaystyle
[X_\ga]\sH|_{e_\gb}=\sum\limits_{S\in\ssyt\spaspb}\,
\prod\limits_{x\in
S}\left(-t_{{\gb'}(x)}-t_{{\gb'}(z(x))}\right)$.
\end{prop}

\begin{ex} Consider $LGr_{3}$, $\ga=\{1,3,\ol{2}\}$, $\gb=\{3,\ol{2},\ol{1}\}$.
Then $\gs(\ga)=(2)$, $\gs(\gb)=(3,2)$, $l(\ga)=2$,
${\gb'}=\{1,2,\ol{3}\}$. The semistandard set-valued tableaux on
$\gs(\gb)$ of shape $\gs(\ga)$ are: \vspace{1em} \psset{unit=1cm}
\begin{equation*}
\pspicture(0,0)(2,1)
\psset{linewidth=.024}
\psline(0,1)(2,1)
\psline(0,0)(2,0)
\psline(0,0)(0,1)
\psline(1,0)(1,1)
\psline(2,0)(2,1)
\rput{0}(.5,.5){$1$}
\rput{0}(1.5,.5){$1$}
\endpspicture
\qquad\qquad\qquad
\pspicture(0,0)(2,1)
\psset{linewidth=.024}
\psline(0,1)(2,1)
\psline(0,0)(2,0)
\psline(0,0)(0,1)
\psline(1,0)(1,1)
\psline(2,0)(2,1)
\rput{0}(.5,.5){$1$}
\rput{0}(1.5,.5){$2$}
\endpspicture
\qquad\qquad\qquad
\pspicture(0,0)(2,1)
\psset{linewidth=.024}
\psline(0,1)(2,1)
\psline(0,0)(2,0)
\psline(0,0)(0,1)
\psline(1,0)(1,1)
\psline(2,0)(2,1)
\rput{0}(.5,.5){$2$}
\rput{0}(1.5,.5){$2$}
\endpspicture
\end{equation*}
\vspace{.2pt}
\begin{equation*}
\pspicture(0,0)(2,1)
\psset{linewidth=.024}
\psline(0,1)(2,1)
\psline(0,0)(2,0)
\psline(0,0)(0,1)
\psline(1,0)(1,1)
\psline(2,0)(2,1)
\rput{0}(.5,.5){$1$}
\rput{0}(1.5,.5){$1,2$}
\endpspicture
\qquad\qquad\qquad
\pspicture(0,0)(2,1)
\psset{linewidth=.024}
\psline(0,1)(2,1)
\psline(0,0)(2,0)
\psline(0,0)(0,1)
\psline(1,0)(1,1)
\psline(2,0)(2,1)
\rput{0}(.5,.5){$1,2$}
\rput{0}(1.5,.5){$2$}
\endpspicture
\end{equation*}
Therefore,
\begin{align*}
[X_\ga]\sK|_{e_\gb}&=
\left(\frac{1}{t_1^2}-1\right)\left(\frac{1}{t_1t_2}-1\right) +
\left(\frac{1}{t_1^2}-1\right)\left(\frac{t_3}{t_2}-1\right) +
\left(\frac{1}{t_2^2}-1\right)\left(\frac{t_3}{t_2}-1\right)\\
&\qquad+
\left(\frac{1}{t_1^2}-1\right)\left(\frac{1}{t_2^2}-1\right)\left(\frac{1}{t_1t_2}-1\right)
+
\left(\frac{1}{t_1^2}-1\right)\left(\frac{1}{t_1t_2}-1\right)\left(\frac{t_3}{t_2}-1\right)
\\
[X_\ga]\sH|_{e_\gb}&=(-2t_1)(-t_1-t_2)+(-2t_1)(-t_2+t_3)+(-2t_2)(-t_2+t_3).
\end{align*}
\end{ex}

\begin{rem}\label{r.positive}
As we shall show in Section \ref{s.class_opp_sch_var}, each term
in the products of Proposition \ref{p.rest_formula_kthry}(i) and
(ii) is of the form  $e^{\gt}-1$ and $\gt$ respectively, where
$\gt$ is a positive root with respect to the Borel subgroup $B^-$.
\end{rem}

\section{The Class of a Schubert Variety}\label{s.class_opp_sch_var}

The \textbf{Pl\"ucker map} $LGr_{n}\to\mathbb{P}(\wedge^n
\bC^{2n})$ is defined by $V\mapsto[v_1\wedge\cdots\wedge v_n]$,
where $\{v_1,\ldots,v_n\}$ is any basis for $V$. The Pl\"ucker map
is a closed immersion, giving $LGr_{n}$ its projective variety
structure.

\subsection*{Reduction to an Affine Variety}

Under the Pl\"ucker map, $e_\gb$ maps to
$[e_{\gb(1)}\wedge\cdots\wedge e_{\gb(n)}]\in\mathbb{P}(\wedge^n
\bC^{2n})$.  Define $p_\gb$ to be the homogeneous
(\textbf{Pl\"ucker}) coordinate $[e_{\gb(1)}\wedge\cdots\wedge
e_{\gb(n)}]^*\in\bC[\mathbb{P}(\wedge^n \bC^{2n})]$.  Let
$\cO_\gb$ be the distinguished open set of $LGr_{n}$ defined by
$p_\gb\neq 0$. Then $\cO_\gb$ is isomorphic to the affine space
$\bC^{n(n+1)/2}$, with $e_\gb$ the origin.  Indeed, $\cO_\gb$ can
be identified with the space of $2n\times n$ complex matrices of
the form $K\cdot M$, where $K$ and $M$ are defined as follows:
\begin{itemize}
\item[1.] $K$ is the $2n\times 2n$ diagonal matrix which has $1$'s
in the diagonal entries of rows $\gb(1),\ldots,\gb(n)$ and rows
$n+1,\ldots,2n$, and $-1$'s in the diagonal entries of all other
rows.

\item[2.] $M$ is any $2n\times n$ complex matrix for which rows
$\gb(1),\ldots,\gb(n)$ form the $n\times n$ identity matrix and
rows $\gb'(1),\ldots,\gb'(n)$ form an $n\times n$ antisymmetric
(i.e., symmetric about the antidiagonal) matrix.
\end{itemize}
%
%
%
Under this identification, we index the rows of $\cO_\gb$ by
$\{1,\ldots,2n\}$ and the columns by $\gb$. We then choose the
coordinates of $\cO_\gb$ to be the matrix elements $y_{ab}$, where
$a\in\gb'$, $b\in\gb$, and $a\leq \ol{b}$ (note: due to the
antisymmetry in the matrices $M$ of 2, each matrix element
$y_{ab}$, where $a\in\gb'$, $b\in\gb$, and $a> \ol{b}$ must be
plus or minus one of our chosen coordinates). Thus
$\{(a,b)\in{\gb'}\times\gb\mid a\leq\ol{b}\,\}$, which we denote
by $\fR_\gb$, forms an indexing set for the coordinates of
$\cO_\gb$.
%
%
\begin{ex}\label{e.big_cell}
Let $n=4$, $\gb=\{1,4,\ol{3},\ol{2}\}$. Then
$\gb'=\{2,3,\ol{4},\ol{1}\}$ and
\begin{equation*}
\cO_\gb=\left\{ \left(
\begin{array}{cccc}
1&0&0&0\\
-y_{21}&-y_{24}&-y_{2\ol{3}}&-y_{2\ol{2}}\\
-y_{31}&-y_{34}&-y_{3\ol{3}}&-y_{2\ol{3}}\\
0&1&0&0\\
y_{\ol{4}\,1}&y_{\ol{4}\,4}&y_{34}&y_{24}\\
0&0&1&0\\
0&0&0&1\\
y_{\ol{1}\,1}&y_{\ol{4}\,1}&y_{31}&y_{21}\\
\end{array}
\right), y_{ab}\in \bC \right\}.
\end{equation*}
\end{ex}
\noindent The space $\cO_\gb$ is $T$-stable, and for ${\bf
s}=\hbox{diag}(s_1,\ldots,s_n,s_{n}^{-1},\ldots,s_{1}^{-1})\in T$
and coordinate functions $y_{ab}\in\bC[\cO_\gb]$,
\begin{equation*}\label{e.char_O_gb}
{\bf s}(y_{ab})=\frac{s_b}{s_a}\,y_{ab},
\end{equation*}
\noindent where $s_{\ol{k}}:=s_k^{-1}$, $k=1,\ldots,n$.

The equivariant embeddings
$e_\gb\stackrel{j}{\to}\cO_\gb\stackrel{k}{\to} LGr_{n}$ induce
homomorphisms
\[
K_T^*(LGr_{n})\stackrel{k^*}{\rightarrow}
K_T^*(\cO_\gb)\stackrel{j^*}{\rightarrow} K_T^*(e_\gb).
\]
The map $j^*$ is an isomorphism, identifying $ K_T^*(\cO_\gb)$
with $K_T^*(e_\gb)$. Define $Y_{\ga,\gb}=X_\ga\cap\cO_\gb$, an
affine subvariety of $\cO_\gb$. We have
\begin{equation*}\label{e.restrict_to_affine}
[X_\ga]\sK|_{e_\gb}=j^*\circ
k^*([X_\ga]\sK)=j^*([k^{-1}X_\ga]\sK)=j^*([Y\ab]\sK)=[Y\ab]\sK.
\end{equation*}
Applying analogous arguments for equivariant cohomology, we obtain
\begin{equation*}\label{e.restrict_to_affine_cohom}
[X_\ga]\sH|_{e_\gb}=[Y\ab]\sH.
\end{equation*}

\subsection*{Reduction to a Union of Coordinate Subspaces}

Let $\spa=\gs(\ga)$, $\spb=\gs(\gb)$. Let $\svt\spaspb$ denote the
set of all set-valued shifted tableaux (not necessarily
semistandard) of shape $\spa$ on $\spb$.  For $S\in\svt\spaspb$,
define
\begin{equation*}\label{e.W_S_defn}
W_S=V(\{y_{{\gb'}(x),\ol{{\gb'}(z(x))}}\mid x\in S\}),
\end{equation*}
a coordinate subspace of $\cO_\gb$.  Define
\begin{equation*}\label{e.W_ab_defn}
W\ab=\bigcup\limits_{P\in\ssyt\spaspb}W_P.
\end{equation*}
The following lemma, whose proof is a consequence of \cite{Gh-Ra}
and appears in Section \ref{s.three_equiv_models}, reduces the
proof of Proposition \ref{p.rest_formula_kthry} to computing the
class of a union of coordinate subspaces.
\begin{lem}\label{l.restrict_to_planes}
$[Y\ab]\sK=[W\ab]\sK$
\end{lem}

\begin{proof}[Proof of Proposition \ref{p.rest_formula_kthry}]
(i) The proofs of Lemmas 4.3 and 4.4 and consequently of
Proposition 2.2(i) of \cite{Kr} carry through if the following
modifications are made: (a) the word `tableau' is replaced by
`shifted tableau' in all steps and all required definitions, and
(b) $t_{\gb(d+1-x)}$ and $t_{{\gb}'(x+c(x)-r(x))}$ are replaced by
$t_{\ol{{\gb'}(z(x))}}$ and $t_{{\gb'}(x)}$ respectively wherever
they occur.  The latter modification accounts for the difference
in the definitions of $W_S$. \vspace{1em}

\noindent (ii) There is a standard ring homomorphism from
$K_T^*(\cO_\gb)$ to $H_T^*(\cO_\gb)$, the Chern character map,
given by $ch: t_i\mapsto e^{-t_i}=1-t_i+t_i^2/2-t_i^3/3+\cdots$.
If $Y\subset \cO_\gb$ is a $T$-stable subvariety, then
\begin{equation*}
ch:[Y]\sK\mapsto [Y]\sH\, +\text{ higher order terms}.
\end{equation*}
Thus, $[Y\ab]\sH$ is the lowest order term of
\begin{equation*}
(-1)^{l(\ga)}\sum\limits_{S\in\ssvt\spaspb}\, \prod\limits_{x\in
S}\left(\frac{1}{e^{-t_{{\gb'}(x)}}e^{-t_{{\gb'}(z(x))}}}-1\right),
\end{equation*}
which equals
\begin{equation*}
\sum\limits_{S\in\ssyt\spaspb}\, \prod\limits_{x\in
S}\left(-t_{{\gb'}(x)}-t_{{\gb'}(z(x))}\right).
\end{equation*}
\end{proof}

\begin{proof}[Proof of Remark \ref{r.positive}]
Let $\pb=\gp(\gb)$, so that
$\spb=\gs(\gb)=\gr(\gp(\gb))=\gr(\pb)$. One can show that
$\pb_j=\#\{i\in\{1,\ldots,n\}\mid {\gb'}(i)< \gb(n+1-j)\}$,
$j=1,\ldots,n$. Therefore
\begin{equation}\label{e.perm_part}
i\leq\pb_j\iff {\gb'}(i)< \gb(n+1-j)
\end{equation}
We look at one term $-t_{{\gb'}(x)}-t_{{\gb'}(z(x))}$ in the
product of (ii). Substituting $i=z(x)$ and $j=x$ into
(\ref{e.perm_part}), we obtain: $z(x)\leq \pb_x\iff
{\gb'}(z(x))<\gb(n+1-x)=\ol{{\gb'}(x)}\iff
{\gb'}(x)<\ol{{\gb'}(z(x))}$. Since $S$ is on $\spb$, $x$
satisfies (\ref{e.semi_stand_on_mu}), i.e., $z(x)\leq
\spb_x+x-1=\pb_x$. Thus ${\gb'}(x)<\ol{{\gb'}(z(x))}$. In
addition, since $x\leq z(x)$, ${\gb'}(x)\leq{\gb'}(z(x))$.

Thus $-t_{{\gb'}(x)}-t_{{\gb'}(z(x))}$ is of the form $-t_a-t_b$,
$a\leq b$, $a<\ol{b}$, and hence $a\leq n$. Clearly this is a
positive root if $b\leq n$.  If $b>n$, then letting $c=\ol{b}$, we
have $-t_a-t_b=-t_a+t_c$, $a< c\leq n$, which is also a positive
root.
\end{proof}

\section{Four Equivalent Models: $\wt{\cF}''\lm$, $\wcF\spaspb$, $\wcD\spaspb$, and $\ssyt\spaspb$}\label{s.three_equiv_models}

In this section we prove Lemma \ref{l.restrict_to_planes}. We
assume all definitions from Sections 5 and 6 of \cite{Kr}. Let
$\pa,\pb$ be symmetric partitions with $\pa\leq\pb$. Let $D_{\pb}$
be the (symmetric) Young diagram associated to $\pb$.
\begin{itemize}
\item A family $F$ of nonintersecting paths on $D_\pb$ is said to
be \textbf{symmetric} if $(i,j)\in F\iff (j,i)\in F$.  In such
case, it can be checked inductively that for any path $p$ of $F$,
$\ol{p}$ is also a path of $F$, where $\ol{p}$ is the path
obtained by replacing each $(i,j)$ of $p$ by $(j,i)$. We define
$\ocF\papb$ and $\ol{\cF''}\!\!\!\papb$ to be the set of all
symmetric elements of $\cF\papb$ and $\cF''\papb$ respectively.

\item A subset $D$ of $D_\pb$ is said to be \textbf{symmetric} if
$(i,j)\in D\iff (j,i)\in D$. We define $\ocD\papb$ to be the set
of all symmetric elements of $\cD\papb$.

\item A semistandard tableau $P$ of shape $\pa$ is said to be
\textbf{symmetric} if $P_{i,j}-i=P_{j,i}-j$ for all $(i,j)\in
D_\pb$. We define $\ossyt\papb$ to be the set of all symmetric
elements of $\text{SSYT}\papb$.
\end{itemize}

\noindent By Lemma 5.14 of \cite{Kr}, $\cF''\papb=\cF\papb$. Hence
$\ol{\cF}\,''\papb=\ol{\cF}\papb$. The bijections
$\cF\papb\to\cD\papb$ and $\cD\papb\to\text{SSYT}\papb$ given in
\cite{Kr} restrict to bijections $\ocF\papb\to\ocD\papb$ and
$\ocD\papb\to\ossyt\papb$ respectively.

Let $\spa=\gr(\pa)$, $\spb=\gr(\pb)$, and let $\wt{D}_{\spb}$ be
the shifted diagram associated with $\spb$. We have the notions of
subsets of $\wt{D}_{\spb}$ and families of nonintersecting paths
on $\wt{D}_{\spb}$, defined analogously as in \cite{Kr}.  If $D\in
\ocD\papb$, then we define $\gr(D)$ to be the subset of
$\wt{D}_{\spb}$ obtained by removing all boxes of $D$ below the
main diagonal of $D$. If $F\in \ocF\papb=\ol{\cF}\,''\papb$, then
we define $\gr(F)$ to be the family of nonintersecting paths on
$\wt{D}_{\spb}$ obtained by removing all boxes in all paths of $F$
below the main diagonal of $F$. If $P\in\ossyt\papb$, then we
define $\gr(P)$ to be the semistandard shifted tableau obtained by
removing all boxes of $P$ below the main diagonal and their
entries. Define $\wcD\spaspb=\gr(\ocD\papb)$,
$\wcF\spaspb=\gr(\ocF\papb)$,
$\wt{\cF}''\spaspb=\gr(\ol{\cF}\,''\papb)$, and note
$\ssyt\spaspb=\gr(\ossyt\papb)$.

We have that $\wt{\cF}''\spaspb=\wt{\cF}\spaspb$, and under $\gr$,
the bijections $\ocF\papb\to\ocD\papb$ and
$\ocD\papb\to\ossyt\papb$ induce bijections
$\wcF\spaspb\to\wcD\spaspb$ and $\wcD\spaspb\to\ssyt\spaspb$
respectively. The following diagram, all of whose squares commute,
summarizes our constructions:
\begin{equation*}
\begin{CD}
\cF''\papb @= \cF\papb @>>> \cD\papb @>>> \text{SSYT}\papb\\
@VVV @VVV @VVV @VVV\\
\ol{\cF}\,''\papb @= \ocF\papb @>>> \ocD\papb @>>> \ossyt\papb\\
@V{\gr}VV @V{\gr}VV @VV{\gr}V @VV{\gr}V\\
\wt{\cF}''\spaspb @= \wcF\spaspb @>>> \wcD\spaspb @>>>
\ssyt\spaspb
\end{CD}
\end{equation*}
All horizontal maps are bijections, as are the four lower vertical
maps. Here we are interested in the bottom row, which gives four
equivalent combinatorial models.

The families $\cF''\papb$ appear in \cite{Ko-Ra}, \cite{Kra1},
\cite{Kra2}, \cite{Kr3}, and \cite{Kr2}; $\cF''\papb$, $\cF\papb$,
$\cD\papb$, and $\text{SSYT}\papb$ appear in \cite{Kr};
$\ol{\cF}\,''\papb$ and $\wt{\cF}''\spaspb$ were introduced in
\cite{Gh-Ra}; $\cD\papb$, $\ocD\papb$, and $\wcD\spaspb$ were
discovered independently by Ikeda-Naruse.

\psset{unit=.29cm}

\begin{ex}\label{ex.dfs} Let $\spa=(3,1)$, $\spb=(5,3,2,1)$.  Below we give all elements
of $\wcF\spaspb$, $\wcD\spaspb$, and $\ssyt\spaspb$.

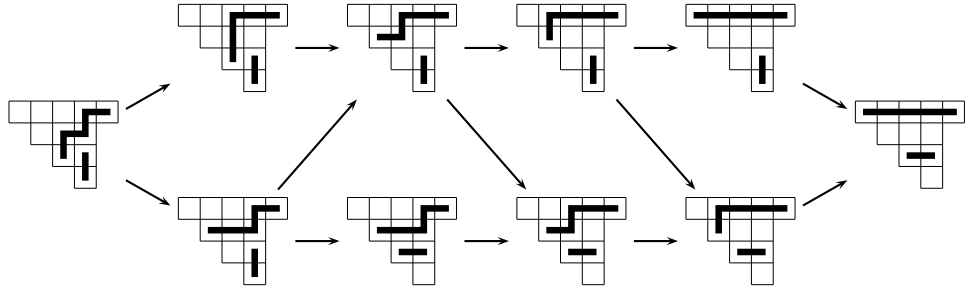
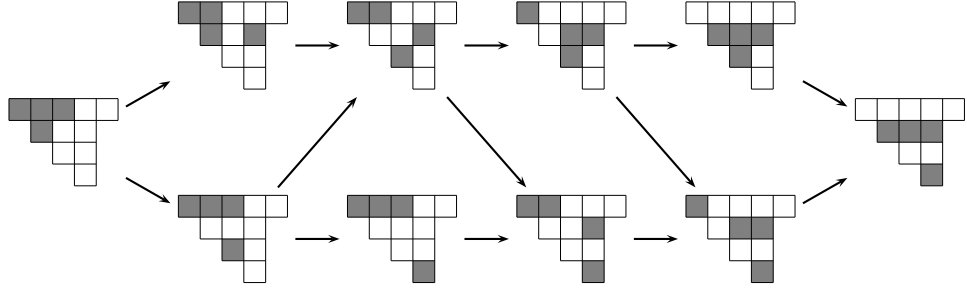
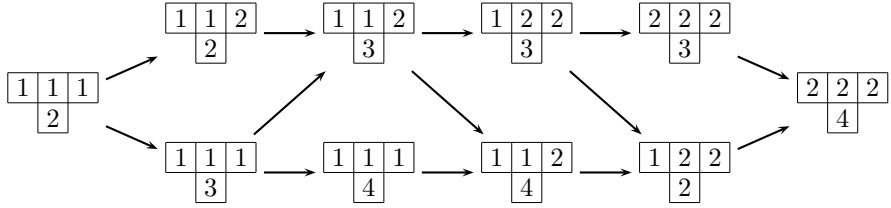
\begin{center}
$ \setlength{\arraycolsep}{.4cm}
\begin{array}{lccccc}
&
\rnode{c}{
\pspicture(0,0)(5,4)
\psset{linewidth=.015}
\psline(0,4)(5,4)
\psline(0,3)(5,3)
\psline(1,2)(4,2)
\psline(2,1)(4,1)
\psline(3,0)(4,0)
\psline(0,3)(0,4)
\psline(1,2)(1,4)
\psline(2,1)(2,4)
\psline(3,0)(3,4)
\psline(4,0)(4,4)
\psline(5,3)(5,4)
\psset{linewidth=.3}
\psline(2.5,1.35)(2.5,3.5)(4.65,3.5)
\psline(3.5,0.35)(3.5,1.65)
\endpspicture
}
&
\rnode{d}{
\pspicture(0,0)(5,4)
\psset{linewidth=.015}
\psline(0,4)(5,4)
\psline(0,3)(5,3)
\psline(1,2)(4,2)
\psline(2,1)(4,1)
\psline(3,0)(4,0)
\psline(0,3)(0,4)
\psline(1,2)(1,4)
\psline(2,1)(2,4)
\psline(3,0)(3,4)
\psline(4,0)(4,4)
\psline(5,3)(5,4)
\psset{linewidth=.3}
\psline(1.35,2.5)(2.5,2.5)(2.5,3.5)(4.65,3.5)
\psline(3.5,0.35)(3.5,1.65)
\endpspicture
}
&
\rnode{e}{
\pspicture(0,0)(5,4)
\psset{linewidth=.015}
\psline(0,4)(5,4)
\psline(0,3)(5,3)
\psline(1,2)(4,2)
\psline(2,1)(4,1)
\psline(3,0)(4,0)
\psline(0,3)(0,4)
\psline(1,2)(1,4)
\psline(2,1)(2,4)
\psline(3,0)(3,4)
\psline(4,0)(4,4)
\psline(5,3)(5,4)
\psset{linewidth=.3}
\psline(1.5,2.35)(1.5,3.5)(4.65,3.5)
\psline(3.5,0.35)(3.5,1.65)
\endpspicture
}
&
\rnode{f}{
\pspicture(0,0)(5,4)
\psset{linewidth=.015}
\psline(0,4)(5,4)
\psline(0,3)(5,3)
\psline(1,2)(4,2)
\psline(2,1)(4,1)
\psline(3,0)(4,0)
\psline(0,3)(0,4)
\psline(1,2)(1,4)
\psline(2,1)(2,4)
\psline(3,0)(3,4)
\psline(4,0)(4,4)
\psline(5,3)(5,4)
\psset{linewidth=.3}
\psline(0.35,3.5)(4.65,3.5)
\psline(3.5,0.35)(3.5,1.65)
\endpspicture
}
&
\\
\rnode{a}{
\pspicture(0,0)(5,4)
\psset{linewidth=.015}
\psline(0,4)(5,4)
\psline(0,3)(5,3)
\psline(1,2)(4,2)
\psline(2,1)(4,1)
\psline(3,0)(4,0)
\psline(0,3)(0,4)
\psline(1,2)(1,4)
\psline(2,1)(2,4)
\psline(3,0)(3,4)
\psline(4,0)(4,4)
\psline(5,3)(5,4)
\psset{linewidth=.3}
\psline(2.5,1.35)(2.5,2.5)(3.5,2.5)(3.5,3.5)(4.65,3.5)
\psline(3.5,0.35)(3.5,1.65)
\endpspicture
}
&
&
&
&
&
\rnode{b}{
\pspicture(0,0)(5,4)
\psset{linewidth=.015}
\psline(0,4)(5,4)
\psline(0,3)(5,3)
\psline(1,2)(4,2)
\psline(2,1)(4,1)
\psline(3,0)(4,0)
\psline(0,3)(0,4)
\psline(1,2)(1,4)
\psline(2,1)(2,4)
\psline(3,0)(3,4)
\psline(4,0)(4,4)
\psline(5,3)(5,4)
\psset{linewidth=.3}
\psline(0.35,3.5)(4.65,3.5)
\psline(2.35,1.5)(3.65,1.5)
\endpspicture
}
\\
&
\rnode{C}{
\pspicture(0,0)(5,4)
\psset{linewidth=.015}
\psline(0,4)(5,4)
\psline(0,3)(5,3)
\psline(1,2)(4,2)
\psline(2,1)(4,1)
\psline(3,0)(4,0)
\psline(0,3)(0,4)
\psline(1,2)(1,4)
\psline(2,1)(2,4)
\psline(3,0)(3,4)
\psline(4,0)(4,4)
\psline(5,3)(5,4)
\psset{linewidth=.3}
\psline(1.35,2.5)(3.5,2.5)(3.5,3.5)(4.65,3.5)
\psline(3.5,0.35)(3.5,1.65)
\endpspicture
}
&
\rnode{D}{
\pspicture(0,0)(5,4)
\psset{linewidth=.015}
\psline(0,4)(5,4)
\psline(0,3)(5,3)
\psline(1,2)(4,2)
\psline(2,1)(4,1)
\psline(3,0)(4,0)
\psline(0,3)(0,4)
\psline(1,2)(1,4)
\psline(2,1)(2,4)
\psline(3,0)(3,4)
\psline(4,0)(4,4)
\psline(5,3)(5,4)
\psset{linewidth=.3}
\psline(1.35,2.5)(3.5,2.5)(3.5,3.5)(4.65,3.5)
\psline(2.35,1.5)(3.65,1.5)
\endpspicture
}
&
\rnode{E}{
\pspicture(0,0)(5,4)
\psset{linewidth=.015}
\psline(0,4)(5,4)
\psline(0,3)(5,3)
\psline(1,2)(4,2)
\psline(2,1)(4,1)
\psline(3,0)(4,0)
\psline(0,3)(0,4)
\psline(1,2)(1,4)
\psline(2,1)(2,4)
\psline(3,0)(3,4)
\psline(4,0)(4,4)
\psline(5,3)(5,4)
\psset{linewidth=.3}
\psline(1.35,2.5)(2.5,2.5)(2.5,3.5)(4.65,3.5)
\psline(2.35,1.5)(3.65,1.5)
\endpspicture
}
&
\rnode{F}{
\pspicture(0,0)(5,4)
\psset{linewidth=.015}
\psline(0,4)(5,4)
\psline(0,3)(5,3)
\psline(1,2)(4,2)
\psline(2,1)(4,1)
\psline(3,0)(4,0)
\psline(0,3)(0,4)
\psline(1,2)(1,4)
\psline(2,1)(2,4)
\psline(3,0)(3,4)
\psline(4,0)(4,4)
\psline(5,3)(5,4)
\psset{linewidth=.3}
\psline(1.5,2.35)(1.5,3.5)(4.65,3.5)
\psline(2.35,1.5)(3.65,1.5)
\endpspicture
}
&
\ncline[nodesep=3pt]{->}{a}{c}
\ncline[nodesep=3pt]{->}{c}{d}
\ncline[nodesep=3pt]{->}{d}{e}
\ncline[nodesep=3pt]{->}{e}{f}
\ncline[nodesep=3pt]{->}{f}{b}
\ncline[nodesep=3pt]{->}{a}{C}
\ncline[nodesep=3pt]{->}{C}{D}
\ncline[nodesep=3pt]{->}{D}{E}
\ncline[nodesep=3pt]{->}{E}{F}
\ncline[nodesep=3pt]{->}{F}{b}
\ncline[nodesep=3pt]{->}{C}{d}
\ncline[nodesep=3pt]{->}{d}{E}
\ncline[nodesep=3pt]{->}{e}{F}
\end{array}
$
\end{center}
\vspace{1em}

\begin{center}
$ \setlength{\arraycolsep}{.4cm}
\begin{array}{lccccc}
&
\rnode{c}{ \pspicture(0,0)(5,4)
\psframe*[linecolor=gray](0,3)(1,4)
\psframe*[linecolor=gray](1,3)(2,4)
\psframe*[linecolor=gray](3,2)(4,3)
\psframe*[linecolor=gray](1,2)(2,3) \psset{linewidth=.03}
\psline(0,4)(5,4) \psline(0,3)(5,3) \psline(1,2)(4,2)
\psline(2,1)(4,1) \psline(3,0)(4,0)
\psline(0,3)(0,4) \psline(1,2)(1,4) \psline(2,1)(2,4)
\psline(3,0)(3,4) \psline(4,0)(4,4) \psline(5,3)(5,4)
\endpspicture
}
&
\rnode{d}{ \pspicture(0,0)(5,4)
\psframe*[linecolor=gray](0,3)(1,4)
\psframe*[linecolor=gray](1,3)(2,4)
\psframe*[linecolor=gray](3,2)(4,3)
\psframe*[linecolor=gray](2,1)(3,2) \psset{linewidth=.03}
\psline(0,4)(5,4) \psline(0,3)(5,3) \psline(1,2)(4,2)
\psline(2,1)(4,1) \psline(3,0)(4,0)
\psline(0,3)(0,4) \psline(1,2)(1,4) \psline(2,1)(2,4)
\psline(3,0)(3,4) \psline(4,0)(4,4) \psline(5,3)(5,4)
\endpspicture
}
&
\rnode{e}{ \pspicture(0,0)(5,4)
\psframe*[linecolor=gray](0,3)(1,4)
\psframe*[linecolor=gray](2,2)(3,3)
\psframe*[linecolor=gray](3,2)(4,3)
\psframe*[linecolor=gray](2,1)(3,2) \psset{linewidth=.03}
\psline(0,4)(5,4) \psline(0,3)(5,3) \psline(1,2)(4,2)
\psline(2,1)(4,1) \psline(3,0)(4,0)
\psline(0,3)(0,4) \psline(1,2)(1,4) \psline(2,1)(2,4)
\psline(3,0)(3,4) \psline(4,0)(4,4) \psline(5,3)(5,4)
\endpspicture
}
&
\rnode{f}{ \pspicture(0,0)(5,4)
\psframe*[linecolor=gray](1,2)(2,3)
\psframe*[linecolor=gray](2,2)(3,3)
\psframe*[linecolor=gray](3,2)(4,3)
\psframe*[linecolor=gray](2,1)(3,2) \psset{linewidth=.03}
\psline(0,4)(5,4) \psline(0,3)(5,3) \psline(1,2)(4,2)
\psline(2,1)(4,1) \psline(3,0)(4,0)
\psline(0,3)(0,4) \psline(1,2)(1,4) \psline(2,1)(2,4)
\psline(3,0)(3,4) \psline(4,0)(4,4) \psline(5,3)(5,4)
\endpspicture
}
&
\\
\rnode{a}{ \pspicture(0,0)(5,4)
\psframe*[linecolor=gray](0,3)(1,4)
\psframe*[linecolor=gray](1,3)(2,4)
\psframe*[linecolor=gray](2,3)(3,4)
\psframe*[linecolor=gray](1,2)(2,3) \psset{linewidth=.03}
\psline(0,4)(5,4) \psline(0,3)(5,3) \psline(1,2)(4,2)
\psline(2,1)(4,1) \psline(3,0)(4,0)
\psline(0,3)(0,4) \psline(1,2)(1,4) \psline(2,1)(2,4)
\psline(3,0)(3,4) \psline(4,0)(4,4) \psline(5,3)(5,4)
\endpspicture
}
& & & & &
\rnode{b}{ \pspicture(0,0)(5,4)
\psframe*[linecolor=gray](1,2)(2,3)
\psframe*[linecolor=gray](2,2)(3,3)
\psframe*[linecolor=gray](3,2)(4,3)
\psframe*[linecolor=gray](3,0)(4,1) \psset{linewidth=.03}
\psline(0,4)(5,4) \psline(0,3)(5,3) \psline(1,2)(4,2)
\psline(2,1)(4,1) \psline(3,0)(4,0)
\psline(0,3)(0,4) \psline(1,2)(1,4) \psline(2,1)(2,4)
\psline(3,0)(3,4) \psline(4,0)(4,4) \psline(5,3)(5,4)
\endpspicture
}
\\
&
\rnode{C}{ \pspicture(0,0)(5,4)
\psframe*[linecolor=gray](0,3)(1,4)
\psframe*[linecolor=gray](1,3)(2,4)
\psframe*[linecolor=gray](2,3)(3,4)
\psframe*[linecolor=gray](2,1)(3,2) \psset{linewidth=.03}
\psline(0,4)(5,4) \psline(0,3)(5,3) \psline(1,2)(4,2)
\psline(2,1)(4,1) \psline(3,0)(4,0)
\psline(0,3)(0,4) \psline(1,2)(1,4) \psline(2,1)(2,4)
\psline(3,0)(3,4) \psline(4,0)(4,4) \psline(5,3)(5,4)
\endpspicture
}
&
\rnode{D}{ \pspicture(0,0)(5,4)
\psframe*[linecolor=gray](0,3)(1,4)
\psframe*[linecolor=gray](1,3)(2,4)
\psframe*[linecolor=gray](2,3)(3,4)
\psframe*[linecolor=gray](3,0)(4,1) \psset{linewidth=.03}
\psline(0,4)(5,4) \psline(0,3)(5,3) \psline(1,2)(4,2)
\psline(2,1)(4,1) \psline(3,0)(4,0)
\psline(0,3)(0,4) \psline(1,2)(1,4) \psline(2,1)(2,4)
\psline(3,0)(3,4) \psline(4,0)(4,4) \psline(5,3)(5,4)
\endpspicture
}
&
\rnode{E}{ \pspicture(0,0)(5,4)
\psframe*[linecolor=gray](0,3)(1,4)
\psframe*[linecolor=gray](1,3)(2,4)
\psframe*[linecolor=gray](3,2)(4,3)
\psframe*[linecolor=gray](3,0)(4,1) \psset{linewidth=.03}
\psline(0,4)(5,4) \psline(0,3)(5,3) \psline(1,2)(4,2)
\psline(2,1)(4,1) \psline(3,0)(4,0)
\psline(0,3)(0,4) \psline(1,2)(1,4) \psline(2,1)(2,4)
\psline(3,0)(3,4) \psline(4,0)(4,4) \psline(5,3)(5,4)
\endpspicture
}
& \rnode{F}{ \pspicture(0,0)(5,4)
\psframe*[linecolor=gray](0,3)(1,4)
\psframe*[linecolor=gray](2,2)(3,3)
\psframe*[linecolor=gray](3,2)(4,3)
\psframe*[linecolor=gray](3,0)(4,1) \psset{linewidth=.03}
\psline(0,4)(5,4) \psline(0,3)(5,3) \psline(1,2)(4,2)
\psline(2,1)(4,1) \psline(3,0)(4,0)
\psline(0,3)(0,4) \psline(1,2)(1,4) \psline(2,1)(2,4)
\psline(3,0)(3,4) \psline(4,0)(4,4) \psline(5,3)(5,4)
\endpspicture
} &
\ncline[nodesep=3pt]{->}{a}{c} \ncline[nodesep=3pt]{->}{c}{d}
\ncline[nodesep=3pt]{->}{d}{e} \ncline[nodesep=3pt]{->}{e}{f}
\ncline[nodesep=3pt]{->}{f}{b}
\ncline[nodesep=3pt]{->}{a}{C} \ncline[nodesep=3pt]{->}{C}{D}
\ncline[nodesep=3pt]{->}{D}{E} \ncline[nodesep=3pt]{->}{E}{F}
\ncline[nodesep=3pt]{->}{F}{b}
\ncline[nodesep=3pt]{->}{C}{d} \ncline[nodesep=3pt]{->}{d}{E}
\ncline[nodesep=3pt]{->}{e}{F}
\end{array}
$
\end{center}
\vspace{1em}

\psset{unit=.4cm}
\begin{center}
$ \setlength{\arraycolsep}{.45cm}
\begin{array}{lccccc}
&
\rnode{c}{\pspicture(0,0)(3,2)
\rput*(.5,1.5){$1$}
\rput*(1.5,1.5){$1$}
\rput*{0}(2.5,1.5){$2$}
\rput{0}(1.5,0.5){$2$}
\psset{linewidth=.03}
\psline(0,2)(3,2)
\psline(0,1)(3,1)
\psline(1,0)(2,0)
\psline(0,1)(0,2)
\psline(1,0)(1,2)
\psline(2,0)(2,2)
\psline(3,1)(3,2)
\endpspicture
}
&
\rnode{d}{
\pspicture(0,0)(3,2)
\rput*(.5,1.5){$1$}
\rput*(1.5,1.5){$1$}
\rput*{0}(2.5,1.5){$2$}
\rput{0}(1.5,0.5){$3$}
\psset{linewidth=.03}
\psline(0,2)(3,2)
\psline(0,1)(3,1)
\psline(1,0)(2,0)
\psline(0,1)(0,2)
\psline(1,0)(1,2)
\psline(2,0)(2,2)
\psline(3,1)(3,2)
\endpspicture
}
&
\rnode{e}{
\pspicture(0,0)(3,2)
\rput*(.5,1.5){$1$}
\rput*(1.5,1.5){$2$}
\rput*{0}(2.5,1.5){$2$}
\rput{0}(1.5,0.5){$3$}
\psset{linewidth=.03}
\psline(0,2)(3,2)
\psline(0,1)(3,1)
\psline(1,0)(2,0)
\psline(0,1)(0,2)
\psline(1,0)(1,2)
\psline(2,0)(2,2)
\psline(3,1)(3,2)
\endpspicture
}
&
\rnode{f}{
\pspicture(0,0)(3,2)
\rput*(.5,1.5){$2$}
\rput*(1.5,1.5){$2$}
\rput*{0}(2.5,1.5){$2$}
\rput{0}(1.5,0.5){$3$}
\psset{linewidth=.03}
\psline(0,2)(3,2)
\psline(0,1)(3,1)
\psline(1,0)(2,0)
\psline(0,1)(0,2)
\psline(1,0)(1,2)
\psline(2,0)(2,2)
\psline(3,1)(3,2)
\endpspicture
}
&
\\
\rnode{a}{
\pspicture(0,0)(3,2)
\rput*(.5,1.5){$1$}
\rput*(1.5,1.5){$1$}
\rput*{0}(2.5,1.5){$1$}
\rput{0}(1.5,0.5){$2$}
\psset{linewidth=.03}
\psline(0,2)(3,2)
\psline(0,1)(3,1)
\psline(1,0)(2,0)
\psline(0,1)(0,2)
\psline(1,0)(1,2)
\psline(2,0)(2,2)
\psline(3,1)(3,2)
\endpspicture
}
&
&
&
&
&
\rnode{b}{
\pspicture(0,0)(3,2)
\rput*(.5,1.5){$2$}
\rput*(1.5,1.5){$2$}
\rput*{0}(2.5,1.5){$2$}
\rput{0}(1.5,0.5){$4$}
\psset{linewidth=.03}
\psline(0,2)(3,2)
\psline(0,1)(3,1)
\psline(1,0)(2,0)
\psline(0,1)(0,2)
\psline(1,0)(1,2)
\psline(2,0)(2,2)
\psline(3,1)(3,2)
\endpspicture
}
\\
&
\rnode{C}{
\pspicture(0,0)(3,2)
\rput*(.5,1.5){$1$}
\rput*(1.5,1.5){$1$}
\rput*{0}(2.5,1.5){$1$}
\rput{0}(1.5,0.5){$3$}
\psset{linewidth=.03}
\psline(0,2)(3,2)
\psline(0,1)(3,1)
\psline(1,0)(2,0)
\psline(0,1)(0,2)
\psline(1,0)(1,2)
\psline(2,0)(2,2)
\psline(3,1)(3,2)
\endpspicture
}
&
\rnode{D}{
\pspicture(0,0)(3,2)
\rput*(.5,1.5){$1$}
\rput*(1.5,1.5){$1$}
\rput*{0}(2.5,1.5){$1$}
\rput{0}(1.5,0.5){$4$}
\psset{linewidth=.03}
\psline(0,2)(3,2)
\psline(0,1)(3,1)
\psline(1,0)(2,0)
\psline(0,1)(0,2)
\psline(1,0)(1,2)
\psline(2,0)(2,2)
\psline(3,1)(3,2)
\endpspicture
}
&
\rnode{E}{
\pspicture(0,0)(3,2)
\rput*(.5,1.5){$1$}
\rput*(1.5,1.5){$1$}
\rput*{0}(2.5,1.5){$2$}
\rput{0}(1.5,0.5){$4$}
\psset{linewidth=.03}
\psline(0,2)(3,2)
\psline(0,1)(3,1)
\psline(1,0)(2,0)
\psline(0,1)(0,2)
\psline(1,0)(1,2)
\psline(2,0)(2,2)
\psline(3,1)(3,2)
\endpspicture
}
&
\rnode{F}{
\pspicture(0,0)(3,2)
\rput*(.5,1.5){$1$}
\rput*(1.5,1.5){$2$}
\rput*{0}(2.5,1.5){$2$}
\rput{0}(1.5,0.5){$2$}
\psset{linewidth=.03}
\psline(0,2)(3,2)
\psline(0,1)(3,1)
\psline(1,0)(2,0)
\psline(0,1)(0,2)
\psline(1,0)(1,2)
\psline(2,0)(2,2)
\psline(3,1)(3,2)
\endpspicture
}
&
\ncline[nodesep=3pt]{->}{a}{c}
\ncline[nodesep=3pt]{->}{c}{d}
\ncline[nodesep=3pt]{->}{d}{e}
\ncline[nodesep=3pt]{->}{e}{f}
\ncline[nodesep=3pt]{->}{f}{b}
\ncline[nodesep=3pt]{->}{a}{C}
\ncline[nodesep=3pt]{->}{C}{D}
\ncline[nodesep=3pt]{->}{D}{E}
\ncline[nodesep=3pt]{->}{E}{F}
\ncline[nodesep=3pt]{->}{F}{b}
\ncline[nodesep=3pt]{->}{C}{d}
\ncline[nodesep=3pt]{->}{d}{E}
\ncline[nodesep=3pt]{->}{e}{F}
\end{array}
$
\end{center}
\end{ex}

\vspace{1em}

\begin{proof}[Proof of Lemma \ref{l.restrict_to_planes}]

Let $\spa=\gs(\ga)$, $\spb=\gs(\gb)$, $\pb=\gp(\gb)$. Recall that
the coordinates $y_{a,b}$ on $\cO_\gb$ are indexed by
$\fR_\gb=\{(a,b)\in{\gb'}\times\gb\mid a\leq\ol{b}\,\}$. Let
$\{v_{a,b}\mid (a,b)\in\fR_\gb\}\subset \cO_\gb$ denote the basis
dual to the basis of linear forms $\{y_{a,b}\mid
(a,b)\in\fR_\gb\}\subset \cO_\gb^*$. For
$F\in\wt{\cF}''_{\spa,\spb}$, define
\begin{equation*}
W_F=\Span(\{v_{{\gb'}(x),\gb(n+1-z)}\mid (x,z)\in
\supp(F)\}\,\dot{\cup}\, \{v_{a,b}\mid (a,b)\in\fR_\gb, a>b\}).
\end{equation*}
In \cite{Gh-Ra}, an explicit equivariant bijection is constructed
from\\ $\bC\left[\bigcup_{F\in\wt{\cF}''\spaspb}W_F\right]$ to
$\bC[Y\ab]$. Thus
\begin{equation}\label{e.ml1}
\Char(\bC[Y\ab])=\Char\left(\bC\left[\bigcup_{F\in\wt{\cF}''\spaspb}W_F\right]\right)
=\Char\left(\bC\left[\bigcup_{F\in\wt{\cF}\spaspb}W_F\right]\right).
\end{equation}

Let $\wt{D}_{\spb}$ be the shifted diagram associated to $\spb$.
For $x,z\in\{1,\ldots,n\}$, $x\leq z$, we have that $(x,z)\in
\wt{D}_{\spb}\iff z\leq\pb_x\iff
{\gb'}(z)<\ol{{\gb'}(x)}\iff{\gb'}(x)<\ol{{\gb'}(z)}=\gb(n+1-z)$
(see proof of Remark \ref{r.positive}). Thus
$\{(a,b)\in\fR_\gb\mid a<b\}$ can be expressed as
$\{({\gb'}(x),\gb(n+1-z))\mid (x,z)\in \wt{D}_{\spb}\}$. Let
$F\in\wcF\spaspb$, and let $D\in \wcD\spaspb$ and $P\in
\ssyt\spaspb$ correspond to $F$ under the bijections above. Since
$\supp(F)$ and $D$ are complements in $\wt{D}_{\spb}$,
\begin{align*} \fR_\gb
&=\{({\gb'}(x),\gb(n+1-z))\mid (x,z)\in \supp(F)\}\,\dot{\cup}\,
\{(a,b)\in\fR_\gb\mid a>b\}\\
&\qquad \dot{\cup}\,\{({\gb'}(x),\gb(n+1-z))\mid (x,z)\in D\}.
\end{align*}
Therefore
\begin{align*}
W_F&=V(\{y_{{\gb'}(x),\gb(n+1-z)}\mid (x,z)\in D\})\\
&=V(\{y_{{\gb'}(x),\gb(n+1-x+r(x)-c(x))}\mid x\in P\})\\
&=V(\{y_{{\gb'}(x),\gb(n+1-z(x))}\mid x\in P\})\\
&=V(\{y_{{\gb'}(x),\ol{{\gb'}(z(x))}}\mid x\in P\})\\
&=W_P.
\end{align*}
Consequently, $\bigcup_{F\in\wt{\cF}\spaspb}W_F =\bigcup_{P\in
\ssyt\spaspb}W_P$, and thus
\begin{equation}\label{e.ml2}
\Char\left(\bC\left[\bigcup_{F\in\wt{\cF}\spaspb}W_F\right]\right)=
\Char\left(\bC\left[\bigcup_{P\in\ssyt\spaspb}W_P\right]\right)=
\Char(\bC[W\ab]).
\end{equation}
Combining (\ref{e.ml1}) and (\ref{e.ml2}), we obtain
$\Char(\bC[Y\ab])=\Char(\bC[W\ab])$. By (4) of \cite{Kr},
$[Y\ab]\sK=[W\ab]\sK$.
\end{proof}


\providecommand{\bysame}{\leavevmode\hbox to3em{\hrulefill}\thinspace}
\providecommand{\MR}{\relax\ifhmode\unskip\space\fi MR }
\providecommand{\MRhref}[2]{%
  \href{http://www.ams.org/mathscinet-getitem?mr=#1}{#2}
}
\providecommand{\href}[2]{#2}

\vspace{1em}

\noindent \textsc{Department of Mathematics, Virginia Tech,
Blacksburg, VA 24063}

\noindent \textsl{Email address}: \texttt{vkreiman@vt.edu}

\vspace{1em}

\noindent July 9, 2006

\end{document}